\documentclass{gtart_h}

\def\ifplaintex{\expandafter\ifx\csname documentclass\endcsname\relax}

\def\gtp{{\mathsurround=0pt\it $\cal G\mskip-2mu$eometry \&\ 
$\cal T\!\!$opology $\cal P\!$ublications}}  % GT publications

\def\recd{{\small Received:\qua\receiveddate\ifx\reviseddate\relax
\else\qquad Revised:\qua\reviseddate\fi\par}} 

\def\lognumber#1{\def\thelognumber{#1}}
\def\volumenumber#1{\def\thevolumenumber{#1}}
\def\volumeyear#1{\def\thevolumeyear{#1}}
\def\papernumber#1{\def\thepapernumber{#1}}
\def\pagenumbers#1#2{\def\startpage{#1}\def\finishpage{#2}}
\def\published#1{\def\publishdate{#1}}

\def\received#1{\def\receiveddate{#1}}
\def\revised#1{\def\reviseddate{#1}}
\def\accepted#1{\def\accepteddate{#1}}

\def\asciiaddress#1{\def\theasciiaddress{#1}}
\def\asciiemail#1{\def\theasciiemail{#1}}
\def\asciiurl#1{\def\theasciiurl{#1}}

\let\\\par\let\thelognumber\relax\let\thevolumenumber\relax
\let\thepapernumber\relax\let\thevolumeyear\relax\let\startpage\relax
\let\finishpage\relax\let\publishdate\relax\let\receiveddate\relax
\let\reviseddate\relax\let\accepteddate\relax\let\theasciititle\relax
\let\theasciiauthors\relax\let\theasciiaddress\relax
\let\theasciiabstract\relax

\let\theasciiemail\relax
\let\theasciiurl\relax

\ifplaintex
\font\logobig=cmssbx10 scaled 3836
\font\logomed=cmssbx10 scaled 2557
\else
\font\logobig=cmssbx10 scaled 4200
\font\logomed=cmssbx10 scaled 2800
\fi

\long\def\makeagttitle{   %%% start of definition of \makeagttitle
\count0=\startpage
\agt\hfill      %   Journal title (top left) 
\hbox to 45truept{\vbox to 0pt{\vglue -13truept{\logomed A\kern -.37em{\logobig 
T}\kern -.38em G}\vss}\hss}
\break
{\small Volume \thevolumenumber\ (\thevolumeyear)
\startpage--\finishpage\nl
Published: \publishdate}

\vglue .25truein

{\parskip=0pt\leftskip 0pt plus
1fil\def\\{\par\smallskip}{\Large\bf\thetitle}\par\medskip} \vglue
0.05truein

{\parskip=0pt\leftskip 0pt plus 1fil\def\\{\par}{\sc\theauthors}
\par\medskip}%
 
\vglue 0.03truein 

{\small\leftskip 25truept\rightskip 25truept{\bf Abstract}\stdspace\theabstract

{\bf AMS Classification}\stdspace\theprimaryclass
\ifx\thesecondaryclass\relax\else; \thesecondaryclass\fi\par
{\bf Keywords}\stdspace \thekeywords\par}\vglue 7truept

}   %%%% end of definition of \makeagttitle

\ifplaintex
\font\phead=cmsl9 scaled 950
\font\pnum=cmbx10 scaled 913
\font\pfoot=cmsl9 scaled 950
\headline{\vbox to 0pt{\vskip -4.5mm\line{\small\phead\ifnum
\count0=\startpage ISSN 1472-2739 (on-line) 1472-2747 (printed)
\hfill {\pnum\folio}\else\ifodd\count0\def\\{ }% 
\ifx\theshorttitle\relax\thetitle\else\theshorttitle\fi\hfill{\pnum\folio}
\else\def\\{ and }{\pnum\folio}\hfill\ifx\theshortauthors\relax\theauthors
\else\theshortauthors\fi\fi\fi}\vss}}
\footline{\vbox to 0pt{\vglue 0mm\line{\small\pfoot\ifnum\count0=\startpage
\copyright\ \gtp\hfill\else
\agt, Volume \thevolumenumber\ (\thevolumeyear)\hfill\fi}\vss}}
\else
\font=cmsl9 scaled 1050
\font\lnum=cmbx10 
\font=cmsl9 scaled 1050
\makeatletter
\def\@oddhead{{\small\lhead\ifnum\count0=\startpage ISSN 1472-2739 
(on-line) 1472-2747 (printed)\hfill {\lnum\number\count0}\else\ifodd\count0
\def\\{ }\ifx\theshorttitle\relax \thetitle \else\theshorttitle\fi\hfill
{\lnum\number\count0}\else\def\\{ and }{\lnum\number\count0}
\hfill\ifx\theshortauthors\relax 
\theauthors\else\theshortauthors\fi\fi\fi}}\def\@evenhead{\@oddhead}
\def\@oddfoot{\small\lfoot\ifnum\count0=\startpage\copyright\ \gtp\hfill\else
\agt, Volume \thevolumenumber\ (\thevolumeyear)\hfill\fi}
\def\@evenfoot{\@oddfoot}
\makeatother
\fi
\let\maketitlepage\makeagttitle
\let\makeshorttitle\maketitlepage
\let\maketitle\maketitlepage

\newwrite\gtoutfile
\long\gdef\makeheadfile{  %%% start of definition of \makeheadfile
{\def\\{, }\def\s{ }
\immediate\openout\gtoutfile head.xxx
\immediate\write\gtoutfile{Proxy-for: \ifx\theasciiauthors\relax
\theauthors\else\theasciiauthors\fi\s<\ifx\theasciiemail\relax\theemail\else\theasciiemail\fi>}
\immediate\write\gtoutfile{\noexpand\\}
\immediate\write\gtoutfile{Authors: \ifx\theasciiauthors\relax
\theauthors\else\theasciiauthors\fi}
{\def\\{ }\immediate\write\gtoutfile{Title: \ifx\theasciititle\relax
\thetitle\else\theasciititle\fi}}
\immediate\write\gtoutfile{Subj-class: GT or SG, GR etc}
\immediate\write\gtoutfile{MSC-class: \theprimaryclass\ifx\thesecondaryclass\relax\else, \thesecondaryclass\fi}
\immediate\write\gtoutfile{Journal-ref: Algebr. Geom. Topol. \thevolumenumber\s
(\thevolumeyear) \startpage-\finishpage}
\immediate\write\gtoutfile{Comments: Published by Algebraic and
Geometric Topology at}
\immediate\write\gtoutfile{\s\s\s  http://www.maths.warwick.ac.uk/agt/AGTVol\thevolumenumber/agt-\thevolumenumber-\thepapernumber.abs.html}
\immediate\write\gtoutfile{\noexpand\\}
\immediate\write\gtoutfile{}
\ifx\theasciiabstract\relax
\immediate\write\gtoutfile{\theabstract}\else
\immediate\write\gtoutfile{\theasciiabstract}\fi
\immediate\write\gtoutfile{}
\immediate\write\gtoutfile{\noexpand\\}
\immediate\write\gtoutfile{}
\immediate\closeout\gtoutfile}}  %%% end of definition of \makeheadfile

\def\maketitlepage{\makeagttitle\makeheadfile}
\let\makeshorttitle\maketitlepage
\let\maketitle\maketitlepage

\lognumber{17}
\volumenumber{4}
\volumeyear{2004}
\papernumber{17}
\published{21 May 2004}
\pagenumbers{311}{332}
\received{17 April 2003}
\revised{26 March 2004}
\accepted{28 April 2004}

\usepackage{amsmath,amssymb}
\usepackage{graphicx}

\newcommand{\kt}{K_t}
\newcommand{\kr}{K_r}
\newcommand{\ym}{\mathcal{YM}}
 
\newcommand{\lcr}{\raisebox{-5pt}{\mbox{}\hspace{1pt}
                 \includegraphics{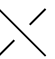}\hspace{1pt}\mbox{}}}
\newcommand{\ift}{\raisebox{-5pt}{\mbox{}\hspace{1pt}
                 \includegraphics{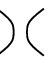}\hspace{1pt}\mbox{}}}
\newcommand{\zer}{\raisebox{-5pt}{\mbox{}\hspace{1pt}
                 \includegraphics{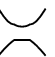}\hspace{1pt}\mbox{}}}
\newcommand{\rkin}{\raisebox{-7pt}{\mbox{}\hspace{1pt}
                 \includegraphics{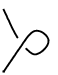}\hspace{1pt}\mbox{}}}
\newcommand{\lkin}{\raisebox{-7pt}{\mbox{}\hspace{1pt}
                 \includegraphics{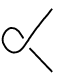}\hspace{1pt}\mbox{}}}
\newcommand{\smst}[1]{\makebox[0pt]{\scriptsize{$#1$}}}
\newcommand{\displaystyl}{\displaystyle}
\newtheorem{theorem}{Theorem}
\newtheorem{lemma}{Lemma}
\newtheorem{cor}{Corollary}
\newtheorem{prop}{Proposition}

\newtheorem{scholium}{Scholium}
\newtheorem{remark}{Remark}

\begin{document}
\title{Shadow world evaluation of the Yang-Mills measure}

\authors{Charles Frohman\\Joanna Kania-Bartoszynska}

\address{Department of Mathematics, University of Iowa, Iowa City, IA
52242, USA\\{\rm and}\\Department of Mathematics, Boise State University, Boise, ID 83725, USA} 
\asciiaddress{Department of Mathematics, University of Iowa, Iowa City, IA
52242, USA\\and\\Department of Mathematics, Boise State University, Boise, ID 83725, USA} 
\gtemail{\mailto{frohman@math.uiowa.edu}\qua {\rm and}\qua
\mailto{kania@math.boisestate.edu}}
\asciiemail{frohman@math.uiowa.edu, kania@math.boisestate.edu}
\gturl{\url{http://www.math.uiowa.edu/~frohman}\qua {\rm and}\qua
\url{http://math.boisestate.edu/~kania}}
\asciiurl{http://www.math.uiowa.edu/ frohman,  http://math.boisestate.edu/ kania} 

\begin{abstract} 
A new state-sum formula for the evaluation of the Yang-Mills measure
in the Kauffman bracket skein algebra  of a closed surface is
derived. The formula extends the Kauffman bracket to diagrams that
lie in surfaces other than the plane.
It also extends Turaev's shadow world invariant of links in a circle bundle
over a surface away from roots of unity.
 The limiting behavior of the Yang-Mills measure when the complex
parameter approaches $-1$ is studied. The formula
is applied to compute integrals
of simple closed curves over  the character variety of the
surface against Goldman's symplectic measure.
\end{abstract}


\maketitle

\section{Introduction}
The Kauffman bracket \cite{Kau} is an unoriented framed version of the Jones
polynomial of oriented links \cite{J}. Kirilov and Reshetikhin gave a state sum
formula for the Kauffman bracket \cite{KR} where the states are admissible
colorings of  the faces
of the link diagram.  This formalism was extended
further to diagrams in surfaces other than the plane by Turaev, who then used
it to give invariants of knots and links in circle bundles over a surface
corresponding to roots of unity \cite{Tu}. This formula is known as {\em  shadow world
evaluation}. 

Let $F$ be an oriented surface and let $t$ be a non-zero complex number.
The Kauffman bracket skein $\kt(F)$ is an algebra derived from
links in a cylinder over $F$. When $t=-1$ this algebra is isomorphic to the 
$SU(2)$-characters of the fundamental group of $F$
\cite{bu},\cite{PS}.
 As $t$ moves away
from $-1$ the deformation of the algebra \cite{quant} corresponds to the Poisson
structure on the characters coming from Goldman's symplectic structure
\cite{Gold1}.

In \cite{nefeli} we constructed a diffeomorphism invariant 
trace $\ym:\kt(F)\rightarrow \mathbb{C}$,
 called the Yang-Mills
measure, that is defined for all
$t$ with $|t|\neq 1$ and when $t$ is a root of unity. 
When $t=-1$ this trace coincides with integration against the symplectic
measure on the  variety underlying the $SU(2)$-characters of
$\pi_1(F)$.  The existence of the trace was established by a study of
its behavior on admissibly colored travalent spines of the surface. 
In this paper we recognize that the Yang-Mills measure can be computed
using a state-sum  formula which corresponds exactly to the formulas
used by Kirilov and Reshetikhin and later by Turaev.

The power of the shadow world formula is that it places the computation of the
Jones polynomial and integration on the $SU(2)$-character variety in the same
framework. Using the combinatorial data from a link diagram in the surface, we
can quickly perform  computations that were before quite 
difficult.  For instance, we prove that
integrating  a simple non-separating curve over the character variety of a
surface always yields zero.
We also  obtain values of the integrals along simple separating curves
on a closed surface of genus bigger than one.

After establishing 
the shadow world formula for the Yang-Mills measure, we use it
together with   recoupling theory  to
compute the limit of the trace at roots of unity 
as the order of the root  tends towards infinity.
When $F$ is a closed surface of genus bigger than $1$ it
converges to twice Goldman's symplectic  measure, i.e., to twice the
classical Yang-Mills measure.

When the complex parameter $t$ is
a root of unity, this formula coincides with Turaev's shadow world
evaluation 
for links in a circle bundle over $F$ with  Euler number
zero \cite{Tu}. 
We  investigate the extent to which Turaev's shadow world evaluation of
links in nontrivial circle bundles can be extended away from roots of unity.
It turns out that the convergence depends on the Euler class of the bundle.

The shadow world evaluation also resembles the formula for the
Turaev-Viro invariants. The estimates in this paper led us to initiate
the study of the Turaev-Viro invariants of three-manifolds away from
roots of unity \cite{normal}. This formula can also be used to
investigate the relationship between the Turaev-Viro invariant and the
Reidemeister torsion of the $SU(2)$-representations of the fundamental
group of the $3$-manifold \cite{F}.

This paper is organized as follows.
Section \ref{kbs} contains definitions and recalls from
\cite{nefeli} the construction of
the Yang-Mills measure away from roots of unity. 
Section \ref{eval} derives the state sum formula for
the Yang-Mills measure.  
In section \ref{shworld} we investigate the extension  of Turaev's
shadow invariant of links in circle bundles over a surface.
In section \ref{lim} we analyze the limit  of the Yang-Mills measure
as the order of roots of unity goes to infinity. 
Finally in the last section  some sample integrals on the 
character variety are evaluated.

This research was partially supported by   NSF grants DMS-0207030 and
DMS-0204627.

\section{ Kauffman bracket skein algebra of a surface}\label{kbs}
\subsection{Preliminaries}
Let $F$ be a compact orientable surface (with or without boundary). 
The cylinder over $F$ is the
three-manifold $F \times [0,1]$.  A framed link  in the cylinder
over  $F$  is an  
embedding of a disjoint union of annuli into $F\times [0,1]$.
Two framed links  are equivalent if there is an 
isotopy of $F\times [0,1]$ taking one to the other.
Framed links in $F \times [0,1]$ can be visualized as diagrams
on the surface $F$. A diagram is an isotopy class of four--valent graphs
in $F$ with a marking at each vertex to indicate which strand goes over
and which strand goes under. The embedded annuli corresponding to
the diagram come from gluing together strips that run parallel to the
surface with the edges of the graph as their core, glued together in
$F \times [0,1]$, according to the  over and under-crossing data.
The  isotopy class of any framed link can be represented this way, and
two diagrams represent the same framed link if and only if they
differ by the second and third Reidemeister moves along with the move
\begin{equation}\rkin=\lkin .\end{equation}
Let $\mathcal{L}$ denote the set of isotopy 
classes of framed links in $F\times [0,1]$, including the empty
link. Fix a complex number $t\neq 0$.  Consider the vector space
$\mathbb{C} \mathcal{L}$ with basis $\mathcal{L}$.  Let $S(F)$ 
be the smallest subspace of $\mathbb{C} \mathcal{L}$ containing all
expressions of the form $\displaystyle{\lcr-t\zer-t^{-1}\ift}$ and
$\bigcirc\cup L+(t^2+t^{-2})L$. The first relation means that there are three framed
links and a ball inside $F\times [0,1]$. The three links coincide outside the
ball, however, the links intersect the ball as blackboard framed arcs
as shown in the relation. The second relation means that there are two framed
links and a ball in the manifold, so that outside the ball the two framed
links coincide, however  one link misses the ball and the other intersects it
in a 0-framed unknot.
The Kauffman bracket
skein  $K_t(F)$ is the quotient
\begin{equation} \mathbb{C} \mathcal{L} / S(F). \end{equation}
The 
algebra structure on $K_t(F)$ comes from laying one link
over the other.  Suppose that $\alpha,\beta \in K_t(F)$ are
skeins represented by links $L_{\alpha}$ and $L_{\beta}$. After
isotopic deformations, to ``lower'' the first link and 
``raise'' the \
second, $L_{\alpha} \subset F\times [0,\frac{1}{2})$, 
and $L_{\beta}
\subset  F\times (\frac{1}{2},1]$. 
The skein $\alpha * \beta $ is
represented by $L_{\alpha} \cup L_{\beta}$. This can be extended using
the distributive law to a 
product on $K_t(F)$.

The notation and the formulas in this paper are taken from
\cite[ch.\ 9]{KL}.
However, the variable $t$ replaces $A$, and we use quantum integers
\begin{equation} [n] = \frac{t^{2n} - t^{-2n}}{t^2-t^{-2}}. \end{equation} When $t=\pm 1$, $[n]=n$.

There is a standard convention, \cite[ch.\ 4]{KL}, for modeling a skein
in $\kt(F)$ on a framed 
trivalent graph  $\Gamma \subset F$.   An
{\em admissible coloring} of $\Gamma$ is an assignment of a
nonnegative integer to each edge of $\Gamma$ so that the colors at trivalent
vertices form admissible triples.
 A triple $(a,b,c)$ is admissible if   $a+b+c$ is even, $a\leq b+c$, $b\leq a+c$,
and $c\leq a+b$. 
  Viewing $F$ as $F\times
\{\frac{1}{2}\}$, the 
skein in $\kt(F)$ corresponding to an admissibly colored graph in $F$  is
obtained by replacing each edge labeled with the 
letter $m$ by the $m$-th Jones--Wenzl
idempotent (see \cite{We}, \cite[ch.\ 3]{KL} or \cite[p.136]{Li}), and
replacing 
trivalent vertices with Kauffman triads (see \cite[ch.\ 4]{KL},
\cite[Fig.\ 14.7]{Li}). 
If $F$ has nonempty boundary and $\gamma$ is a trivalent spine of $F$, then
the collection of skeins modeled on all admissible colorings of $\gamma$ forms
a basis for $\kt(F)$.

\begin{figure}[ht!]
\begin{center}
\begin{picture}(56,120)
\raisebox{8pt}{\hspace{-1in}\includegraphics{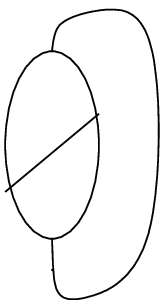}
\put(-50,16){$a$}
\put(-60,49){$b$}
\put(5,49){$e$}
\put(-23,70){$c$}
\put(-15,40){$d$}
\put(-40,49){$f$}}
\end{picture}
\begin{picture}(46,45)
\hspace{0.1in}\raisebox{40pt}{\includegraphics{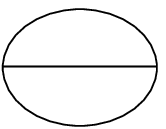}
\put(-25,5){$a$}
\put(-25,19){$b$}
\put(-25,37){$c$}}
\end{picture}
\end{center}
\caption{Tet and theta}\label{tet}
\end{figure}

We often use values of the Kauffman bracket of the following three
skeins in a disk.
The zero framed unknot colored with the
$n$-th Jones-Wenzl  
idempotent, denoted by $\Delta_n$, is
equal to $(-1)^n[n+1]$. 
The symbol $\text{Tet}\begin{pmatrix} a & b & e \\ c & d &
f\end{pmatrix}$ stands for the bracket of the skein modeled by the
colored graph  pictured 
in Figure \ref{tet} on the left.  Finally, the quantity
$\theta(a,b,c)$ is the bracket 
of the skein pictured on the right in Figure
\ref{tet}. 
The explicit formulas are given in
\cite[pp.97-8]{KL}. 

Recall the fusion identity:
\begin{equation}
\raisebox{-24pt}{\includegraphics{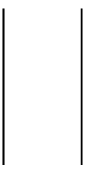}}
\hspace{-20pt}
%\raisebox{27pt}{\smst{}\hspace{27pt}\smst{b}}
\raisebox{23pt}{\smst{a}\hspace{23pt}\smst{b}}
= 
\sum_{c}  \frac{\Delta_c}{\theta(a,b,c)} 
\raisebox{-24pt}{\includegraphics{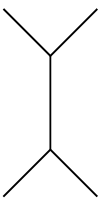}}
\hspace{-9pt} 
\raisebox{3pt}{\smst{c}} 
\hspace{-9pt}
\raisebox{27pt}{\smst{a}\hspace{22pt}\smst{b}}
\hspace{-22pt}
\raisebox{-22pt}{\smst{a}\hspace{22pt}\smst{b}}
\end{equation}
where the sum is over all $c$ so that the triples $(a,b,c)$ are
admissible.
The fusion relation is satisfied in $\kt(F)$ unless
$t$ is a root of unity other than $\pm 1$.

When $t=e^{\frac{\pi i}{2r}}$ for  odd $r>1$ we work in the quotient
of $\kt(F)$ where the appropriate form of the fusion identity has been
enforced. In this case triples are 
$r$-admissible, that is $a+b+c\leq2r-4$, in addition to the 
previous conditions. Dividing by the fusion relation is equivalent to
taking a quotient of $\kt(F)$ by the submodule spanned by skeins
corresponding to trivalent graphs where some edge carries the label $r-1$
\cite{Ro}.
We denote the quotient algebra
$\kr(F)$. Just as admissibly colored spines of $F$ form a basis for $\kt(F)$,
$r$-admissibly colored spines are a basis for $\kr(F)$. 
If $F$ is closed, $\kr(F)$ is isomorphic to $M_n(\mathbb{C})$ 
for some $n$ that
depends on $r$ and on the topological type of the surface \cite{Ro}. 

The $6j$-symbol is defined by
\begin{equation} \left\{ \begin{matrix} a & b & e \\ c & d & f\end{matrix}\right\} = \frac{\text{Tet}\begin{pmatrix} a & b & e \\ c & d &
f\end{pmatrix}\Delta_e}{\theta(a,d,e)\theta(b,c,e)}.\end{equation}
The recoupling formula says that:

\vspace{.13in}

\hspace{2in}\raisebox{-10pt}{\begin{picture}(25,25)\includegraphics{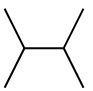} 
\put(0,-4){$d$}\put(0,25){$c$}\put(-29,-4){$a$}\put(-29,25){$b$}
\put(-14,15){$f$}\end{picture}} \  
$= \displaystyle\sum_e \left\{ \begin{matrix} a & b & e \\ c
    & d & f\end{matrix}\right\}$ \raisebox{-10pt}{\begin{picture}(25,25) 
\includegraphics{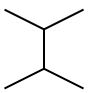}
\put(0,-4){$d$}\put(0,25){$c$}\put(-29,-4){$a$}\put(-29,25){$b$}
\put(-10,11){$e$}\end{picture}}.

\vspace{.13in}

All of the quantities $\Delta_a$, $\theta(a,b,c)$,
\[\text{Tet}\begin{pmatrix} a & b & e \\ c & d &
f\end{pmatrix},\ \text{and}
 \left\{ \begin{matrix} a & b & e \\ c & d & f\end{matrix}\right\}\] may
fail to be defined at $4r$th roots of unity when one of the arguments
is bigger than $r-2$. However, fixing the arguments, and viewing them
as functions of $t$, they are all continuous in neighborhoods of
$t=-1$. 

Given  
a link $L\subset F\times I$ and a spine $\gamma$ of $F$, 
for each $t$ (and for each $r$) $L$ can be written as a linear combination of
skeins corresponding to admissible (and respectively $r$-admissible)
colorings of the  spine.  These coefficient are not necessarily
continuous as $t$ approaches $e^{\frac{\pi i}{2r}}$. However, they are
continuous at $t= -1$.

\subsection{Yang-Mills measure} 
Suppose that $F$ is a closed surface of genus $g$. 
Since the quotient  of $M_n(\mathbb{C})$ by the subvector
space of commutants is one dimensional
there is a unique trace
$\ym:\kr(F)\rightarrow \mathbb{C}$
so that $\ym(\emptyset)=\sum_{i=0}^{r-2} \frac{1}{[i+1]^{(2g-2)}}$. 
By the trace we mean that $\ym(\alpha* \beta)=\ym(\beta*\alpha)$ for all
$\alpha,\beta$.
This trace is called the Yang Mills measure. If $t=\pm 1$ or
$|t|\neq 1$, there is a  diffeomorphism invariant,
locally defined trace
on $K_t(F)$ so that when $F$ is closed and $g>1$,
$\ym(\emptyset)=\sum_{i=0}^{\infty} \frac{1}{[i+1]^{(2g-2)}}$, and 
when $g=1$ or when $F$ has boundary then $\ym(\emptyset)=1$, \cite{nefeli}. 
If the surface $F$ has boundary, the locality condition can be characterized in
the following way.
If $\kappa$ is a properly embedded
arc in $F$, let $F'$ be the result of cutting $F$ along $\kappa$. Let
$N_t(F)$ be the linear subspace of $\kt(F)$ spanned by skeins that
are represented by admissibly colored trivalent graphs that intersect
$\kappa$ in a single point of transverse intersection, so that the
edge intersecting $\kappa$ carries a nonzero label. Let
$\iota :\kt(F')\rightarrow \kt(F)$ be the natural inclusion map.
Then 
\begin{equation}\label{iota}
\kt(F)=\mathrm{im}(\iota) \oplus N_t(F).
\end{equation}
{\em Locality} is the statement
that $\ym$ is zero on $N_t(F)$ and its restriction to $\mathrm{im}(\iota)$
is equal to $\ym$ on $\kt(F')$.

Let $A$ be an annulus which is the result of removing a regular
neighborhood of a point from the
disk $D$, and let 
$s_u$ be the skein in $A$ obtained by coloring a core of the annulus
with the $u$-th Jones-Wenzl idempotent. Given a skein $\alpha$ in $D$, perturb
the links representing it so that they miss the removed disk to
get a skein in $A$, denoted $\alpha'$ .
For the remainder of this and in the next section, use subscripts
to indicate the domain of $\ym$.  
It is proved in \cite{nefeli} that 
\begin{equation}\label{annulus}
\ym_D(\alpha)=\sum_{u=0}^{\infty} \Delta_u \ym_A(\alpha'*s_u).
\end{equation}
Although the sum is infinite only a finite number of terms are nonzero.

If $F$ is a closed surface of genus greater than 1,
let $F'$ be the result of removing a regular neighborhood of a  point
from $F$. Let $\partial_u$ be a small zero framed knot parallel to the
boundary of $F'$, colored
with the $u$-th Jones-Wenzl idempotent. 
For a skein $\alpha$  in $F$,
let $\alpha'$ be a skein in $F'$ obtained by perturbing all the links in
a representative of $\alpha$ to lie in $F'$.
Then
\begin{equation}\label{defym} 
\ym_F(\alpha)=\sum_{u=0}^{\infty} \Delta_u \ym_{F'}(\alpha'*\partial_u).
\end{equation}
The series is absolutely convergent \cite{nefeli} if
$|t|\neq 1$ and when $t=\pm 1$. 

When $F$ has genus
 $g=1$ we define the quantities as above except that
the Yang-Mills measure is given by the formula,
\begin{equation} \ym_F(\alpha)=\lim_{n\rightarrow \infty} \frac{1}{n}\sum_{u=0}^n
\Delta_u \ym_{F'}(\alpha'*\partial_u).\end{equation} 
Recall that the Kauffman bracket skein module of a disjoint union of
two $3$-manifolds is the tensor product of the modules of the
pieces. In particular if a surface $F$ is a disjoint union of surfaces
$F_1$ and $F_2$ then
\begin{equation}
\kt(F)=\kt(F_1)\otimes \kt(F_2).
\end{equation}
The Yang Mills measure is bilinear with respect to the disjoint
union. That is, if $\alpha\in\kt(F)$, where $F=F_1\sqcup F_2$ and
$\alpha=\alpha_1*\alpha_2$ 
whith $\alpha_1\in\kt(F_1)$ and $\alpha_2\in\kt(F_2)$ then  
$\ym_t(\alpha)=\ym_t(\alpha_1)\ym_t(\alpha_2)$.

\section{A formula for evaluating the Yang-Mills measure of a colored
  framed link  }\label{eval}

\begin{figure}[ht!]
\begin{center}
\begin{picture}(31,31)
\includegraphics{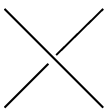}
\put(-21,-2){$\frac{1}{2}$}
\put(-21,27){$\frac{1}{2}$}
\put(-10,13){$-\frac{1}{2}$}
\put(-40,13){$-\frac{1}{2}$}
\end{picture}
\end{center}
\caption{Contributions to the gleam at a vertex}\label{versigns}
\end{figure}

We begin with notation for describing the geography of link
diagrams. The projection of a link onto the surface decomposes it 
into faces, edges
and vertices.  Given a face $f$, at each of its vertices there are two
edges lying on the boundary of $f$; one of them is an over, and the
other an under-crossing. Associate $-\frac{1}{2}$ to the vertex when in
the counter-clockwise order the under comes before the over-crossing
edge. Otherwise assign $\frac{1}{2}$ to the vertex (see Figure
\ref{versigns}). The gleam $g_f$ of the face is the sum of these
assignments over all the vertices lying on the boundary of $f$.
The Euler characteristic of the face $f$ is  denoted
$x_f$. 

Each edge $e$ lies between two faces $f_{e,1}$ and $f_{e,2}$.
However, not every edge must belong to a vertex. Let $v(e)=1$ if $e$
is adjacent to some vertex, and $v(e)=0$ otherwise.
The edge $e$ carries a color $k_e$ inherited from the coloring of the link.

\begin{figure}[ht!]
\begin{center}
\begin{picture}(31,31)
\includegraphics{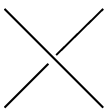}
\put(-25,0){$f_{v,1}$}
\put(-25,27){$f_{v,3}$}
\put(-6,13){$f_{v,4}$}
\put(-38,13){$f_{v,2}$}
\end{picture}
\end{center}
\caption{Labeling faces around a vertex}\label{verfaces}
\end{figure}

At each vertex $v$ there are four edges 
corresponding to at most two components of the link, the over-crossing one colored
with $k_{v,o}$, and the under-crossing one colored with $k_{v,u}$.
Any vertex $v$ is adjacent to four corners of faces, labeled
clockwise
$f_{v,1}$, $f_{v,2}$, 
$f_{v,3}$ and $f_{v,4}$, 
so that the under-crossing edge lies between $f_{v,1}$ and $f_{v,2}$
as pictured in Figure \ref{verfaces}.

\begin{figure}[ht!]
\begin{center}
\begin{picture}(55,75)
 \scalebox{.5}{\includegraphics{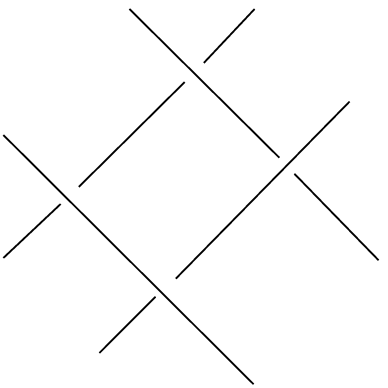}}
\put(0,12){$1$}
\put(0,37){$3$}
\put(-15,50){$2$}
\put(-15,0){$4$}
\put(-28,52){$v$}
\put(-45,40){$f_{1}$}
\put(-32,27){$f_{4}$}
\end{picture}
\end{center}
\caption{Fragment of a colored link diagram}\label{diag}
\end{figure}

An admissible coloring of the diagram is an assignment of
a color $u_f \in \mathbb{N}$ to each face $f$ so that along each
edge the colors $u_{f_{e,1}}$,$u_{f_{e,2}}$, and $k_e$ form an
admissible triple. 
Figure \ref{diag} shows a piece of a colored link diagram. The edges
carry colors $1$, $2$, $3$ and $4$. The vertex $v$ has $k_{v,o}=1$ and
 $k_{v,u}=2$. Note that admissibility forces $u_{f_1}$ to be equal 
$u_{f_4}-2$, $u_{f_4}$ or $u_{f_4}+2$. The gleam of $f_4$ is $1$. 

\begin{theorem}\label{theformula} 
Let $F$ be a closed surface of genus $g>1$, and $\mathcal{L}$ 
 a framed colored link in $ F \times I$, whose diagram has faces
$f$, edges $e$ and vertices $v$.  Suppose that
$|t|\neq 1$ or $t=\pm 1$. The value of $\ym:\kt(M)\rightarrow \mathbb{C}$ on
the skein corresponding to $\mathcal{L}$ is given by
\begin{equation}\label{shadow}
\sum_{u_f}(-1)^{\sum_fg_fu_f}t^{g_fu_f(u_f+2)}\frac{
\prod_f \Delta_{u_f}^{x_f}\prod_v\mathrm{Tet}
\begin{pmatrix} u_{f_{v,1}} & u_{f_{v,2}} &
  k_{v,u} \\ u_{f_{v,3}} & 
  u_{f_{v,4}} & k_{v,o}\end{pmatrix}}
{\prod_e\theta(u_{f_{e,1}},u_{f_{e,2}},k_e)^{v(e)}},
\end{equation}
where the sum is over all admissible colorings $u_f$  of the
faces.
\end{theorem}

\proof
The computation of the Yang-Mills
measure of a skein  on a closed surface starts by removing a 
disk from the surface
and then taking the sum of values of the Yang-Mills measure of an
associated family of
skeins on the punctured surface. The first step of the proof is
to  show  
that an arbitrary number of punctures can be used.
In this way we  obtain a multiply
infinite sum whose value is the Yang-Mills measure. The second step is
reducing this sum to the expression (\ref{shadow}).

We prove the first step in a sequence of lemmas.
\begin{lemma}\label{inadisk}
Let $\{p_i\}_{i=1}^n$ be a collection of points in the interior of a
disk $D$. Denote by $D'$ the planar surface obtained by removing
disjoint regular neighborhoods of the points $p_i$ from $D$, and let
$L\subset D\times [0,1]$ be a framed link whose diagram lies in $D'$.
For each $i=1,\dots,n$, let $\partial_{i,u}$ be a blackboard framed knot
parallel to the boundary of the neighborhood of $p_i$, colored with
the $u$-th Jones-Wenzl idempotent. Then
\begin{equation}\ym_D(L)=\sum_{u_1=0}^{\infty}\dots\sum_{u_n=0}^{\infty}
(\prod_{i=1}^n\Delta_{u_i})
\ym_{D'}(L*\prod_{i=1}^n\partial_{i,u_i}).\end{equation}
\end{lemma}
\proof By induction on $n$.
Even though the sum is infinite, only finitely many terms are nonzero.
For $n=1$ this is
identity (\ref{annulus}). 
Suppose the lemma is true for $n-1$ points, and choose a proper arc
$\kappa\subset D'$ that separates point $p_n$ from
$\{p_1,\dots,p_{n-1}\}$. Let $D_1$ and $D_2$ be the components of
$D\setminus\kappa$ with $\{p_1,\dots,p_{n-1}\}\subset D_1$, and let
$D_1'$, $D_2'$ be the corresponding components of $D'\setminus\kappa$.
Recall that by (\ref{iota}), the  skein $[L]$ corresponding to the link $L$
can be written as $\lambda_0+\lambda_n$, where 
$\lambda_0\in \mathrm{im}\left(\iota:K_t(D_1)\otimes K_t(D_2)
  \rightarrow K_t(D)\right)$, and $\lambda_n\in N_t(D)$. We can
express $\lambda_0$ as
\begin{equation}\lambda_0=\iota\left(\sum_j\lambda_{1,j}\otimes \lambda_{2,j}\right).\end{equation}
Thus by locality and by linearity,
\begin{multline}\ym_D(L)= \\
\sum_j\sum_{u_1=0}^{\infty}\dots\sum_{u_n=0}^{\infty}
(\prod_{i=1}^n\Delta_{u_i})
\ym_{D_1}(\lambda_{1,j}*\prod_{i=1}^{n-1}\partial_{i,u_i})
\ym_{D_2}\left(\lambda_{2,j}*\partial_{n,u_n}\right).\end{multline}
By the distributive law this is equal to
\begin{multline}\sum_j 
\left(\sum_{u_1=0}^{\infty}\dots\sum_{u_{n-1}=0}^{\infty}
\prod_{i=1}^{n-1}\Delta_{u_i}
\ym_{D'_1}(\lambda_{1,j}*\prod_{i=1}^{n-1}\partial_{i,u_i})\right)
\cdot\\
\left(\sum_{u_n=0}^{\infty}\Delta_{u_n}
\ym_{D'_2}(\lambda_{2,j}*\partial_{n,u_n})  \right).
\end{multline}
By induction this is further equal to
\begin{equation}\sum_j\ym_{D_1}(\lambda_{1,j})\ym_{D_2}(\lambda_{2,j}),\end{equation}
which, by locality, is $\ym_D(L)$.\endproof

Next, we generalize Lemma \ref{inadisk} to an arbitrary surface with
boundary. 
\begin{lemma}\label{withboundary}
Let  $\{p_i\}_{i=1}^n$ be a collection of points in the interior of a
compact oriented surface $F$ with boundary.
Denote by $F'$ the result of removing
disjoint regular neighborhoods of the points $p_i$ from $F$, and let
$L\subset F\times [0,1]$ be a framed link whose diagram lies in $F'$.
The collection of knots $\partial_{i,u_i}$ is defined as in Lemma
\ref{inadisk}. Then
\begin{equation}\ym_F(L)=\sum_{u_1=0}^{\infty}\dots\sum_{u_n=0}^{\infty}
(\prod_{i=1}^n\Delta_{u_i})
\ym_{F'}(L*\prod_{i=1}^n\partial_{i,u_i}).\end{equation}
\end{lemma}
\proof
Choose a family of proper arcs (called cross-cuts) in $F$
which lie in $F'$ and 
cut $F$ into a disk $D$. Using identity (\ref{iota}) to decompose the
skein $[L]$, with respect to the crosscuts we get
 $[L]=\lambda_0+\lambda_n$  where
\begin{equation}\ym_{F'}(L)=\ym_D(\lambda_0).\end{equation}
Now apply Lemma \ref{inadisk}.\endproof

Finally, we prove the analogous lemma for closed surfaces.

\begin{lemma}
Suppose that $F$ is a closed surface of genus $g>1$. Let 
$\{p_i\}_{i=1}^n$, $F'$, $L$,  and 
$\partial_{i,u_i}$ be  as in Lemma \ref{withboundary}. Then 
\begin{equation}\ym_F(L)=
\sum_{u_1=0}^{\infty}\dots\sum_{u_n=0}^{\infty}
(\prod_{i=1}^n \Delta_{u_i})
\ym_{F'}(L*\prod_{i=1}^n\partial_{i,u_i}).\end{equation}
\end{lemma}
\proof
Let $p_0$ be a point in $F'$ disjoint from the diagram of $L$, and $N(p_0)$ its regular neighborhood. By the 
definition of the Yang-Mills measure (\ref{defym}),
\begin{equation}\ym_F(L)=
\sum_{u_0=0}^{\infty}\Delta_{u_0}
\ym_{F\setminus N(p_0)}(L*\partial_{0,u_0}).\end{equation}
Now apply Lemma \ref{withboundary}.\endproof

Thus the value of the Yang-Mills measure of a link can be computed
using an arbitrary number of punctures.
In order to prove Theorem \ref{theformula}, choose one point $p_i$ in
each face $f_i$ of the diagram of the link $L$.
The next step is to apply locality in the multiply punctured surface $F'$ to
compute the measure.

\begin{figure}[ht!]
\[ \ym\left( \raisebox{-10pt}{\includegraphics{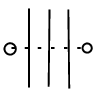}}\right)
=\sum_c \frac{\Delta_c}{\theta(c,u_{i},k_e)}
\ym\left( \raisebox{-10pt}{\includegraphics{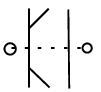}}\right)\]
\[
=\sum_{c,d}\frac{\Delta_c \Delta_d}{\theta(c,u_{i},k_e) \theta(c,d,u_{j})}
\ym\left(
  \raisebox{-10pt}{\includegraphics{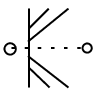}}\right)
=\frac{1}{\theta(u_{i},k_e,u_{j})}\ym\left(
  \raisebox{-10pt}{\includegraphics{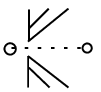}}\right).\]
\caption{Fusing along cross-cuts.}\label{crosscuts}
\end{figure}

Choose a system of cross-cuts joining  the punctures of adjacent faces,
along with $2h+b-1$
cross-cuts in the interior of each face of genus $h$ and $b$ boundary
components, to cut it down to a disk.
Fuse in neighborhoods of cross-cuts (as in Figure \ref{crosscuts}) to
write the skein $[L]$ as a linear 
combination of colored trivalent graphs, each intersecting every
cross-cut at most once transversally. 
In  Figure \ref{crosscuts} the crosscuts are depicted by dotted
horizontal  lines joining the little circles (i.e. boundaries of the
punctures). The middle vertical line is a strand of the link, colored
with $k_e$, the left vertical line is a strand of the knot
$\partial_i$ colored with $u_i$, and the right vertical line is  a
strand of the knot $\partial_j$ colored with $u_j$.

\begin{figure}[ht!]
\begin{center}
\begin{picture}(85,70)
\scalebox{1.75}{\includegraphics{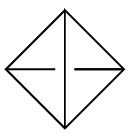}}
\put(-70,52){$u_{f_{v,2}}$}
\put(-10,10){$u_{f_{v,4}}$}
\put(-10,52){$u_{f_{v,3}}$}
\put(-70,10){$u_{f_{v,1}}$}
\put(-50,37){$k_{v,o}$}
\put(-30,23){$k_{v,u}$}
\end{picture}
\end{center}
\caption{Tetrahedron around a vertex}\label{tetr} 
\end{figure}

The measure of $[L]$ is equal to the sum of the 
measures of those terms where every edge intersecting a cross-cut has
color zero, with additional coefficients coming from fusion.
These terms are unions of tetrahedra and theta curves. 
The tetrahedra correspond to the crossings of $L$, i.e., to the
vertices of the diagram of $L$  (see
figure \ref{tetr}). The
theta curves come from strands of the link which are not involved in
any crossings, i.e., from edges of the link diagram which are not
adjacent to any vertices. 

Fusing along each of the $2h+b-1$ cross-cuts 
for each face $f$ introduces the
coefficient $\frac{1}{\theta(u_{f},
  u_{f},0)}=\frac{1}{\Delta_{u_{f}}}$.
Recalling the coefficient $\Delta_{u_{f}}$ coming from
puncturing the face $f$ and totaling this with the contribution from
the cross-cuts yields
that the total factor coming from the face $f$ is
$\left(\Delta_{u_{f}}\right)^{x_{f}}$.

Fusing along a cross-cut joining two faces $f_1$, $f_2$ that intersect
along an edge $e$ introduces the coefficient
$\frac{1}{\theta(u_{f_{e,1}},u_{f_{e,2}},k_e)}$. If $e$ is not adjacent to any
vertex this coefficient cancels with the measure of the resulting
theta curve.

To compute the measure of a tetrahedron at the vertex $v$ slide
the top edge in figure \ref{tetr} to the side and use the tulip identity
\cite[Fig.\ 14.14]{Li}). This evaluates to
\begin{equation} (-1)^{n_1} 
t^{n_2}  
 \mathrm{Tet}\begin{pmatrix} u_{f_{v,1}} & u_{f_{v,2}}
  & k_{v,u} \\ u_{f_{v,3}} & u_{f_{v,4}} & k_{v,o}\end{pmatrix},\end{equation}  
where $n_1=\frac{u_{f_{v,1}}+u_{f_{v,3}}-u_{f_{v,2}}-u_{f_{v,4}}}{2}$,
and \[ n_2=\frac{-u_{f_{v,1}}(u_{f_{v,1}}+2)-u_{f_{v,3}}(u_{f_{v,3}}+2)+
    u_{f_{v,2}}(u_{f_{v,2}}+2)+u_{f_{v,4}}(u_{f_{v,4}}+2)} 
  {2}.\]
Taking the product of all these,  collecting the exponents  of $t$
according
to the faces, and collecting the powers  of $-1$ according to the vertices,
yields the  desired formula.
\endproof

\begin{scholium}
The analogous formula holds when $t=e^{\frac{\pi i}{2r}}$ for odd 
$r\geq 3$, with the altered requirement  that the colors $u_f$ be
$r$-admissible. \qed
\end{scholium}

Let $F$ be a compact oriented surface with boundary and $L$ a blackboard 
framed colored link in $F \times [0,1]$. An admissible coloring of a diagram 
of $L$ is an assignment of colors to the faces of the diagram so that each 
face 
containing a boundary component of $F$ carries the label $0$ and along each 
edge the triple of colors coming from the edge and the adjacent faces is 
admissible.
\begin{theorem} If $ F$  is a compact oriented surface with boundary then
the Yang-Mills measure of a skein corresponding to a blackboard framed
colored link can be computed using formula 
(\ref{shadow}) where the sum is over all admissible colorings of a
diagram  corresponding to the link. \qed
\end{theorem}

\section{The Shadow World}\label{shworld}

Turaev \cite{Tu} constructed a formalism for invariants of framed links
in a circle bundle over a surface $F$. 
A {\em shadow} $S$ is a four-valent
graph in a closed, oriented surface, with an assignment of an integer
to each face 
of the graph,  called the gleam of the face. The sum of
all the gleams over all the faces is the total gleam of the shadow.
Let $p:E\rightarrow F$ be a circle bundle. If $L\subset E$ is a link in
general position with respect to $p$, then its image in $F$ is a four-valent
graph. There is a recipe for deriving gleams from how the link lies in
$E$. The total gleam
is the negative of the Euler number of the circle bundle that 
the framed link
lives in. 

One can assign a shadow to a framed link in $F\times I$. Make sure
it has  the blackboard framing and then project it down to the surface.
The graph is the image of the link under
projection.  The gleam of each face is assigned as in
the computation of the
Yang-Mills measure (see Figure \ref{versigns}). 
This coincides with the rule given in the Remark on
page 41 of \cite{Tu}.

If the link is colored, the edges of its shadow inherit the colors.
An admissible coloring of a shadow $S$ is an assignment of colors to the
faces so that for each edge of the shadow the triple of colors from
the adjacent faces and the edge is admissible. The formula for
Turaev's shadow world evaluation of $S$ is a sum analogous to
(\ref{shadow}). However, there is no notion of over-crossing and
under-crossing at a vertex. Hence for the gleams we use the assigned
values $g_f$. Suppose that the colors of faces at a vertex are  $u_{f_{v,1}}$,
$u_{f_{v,2}}$,  $u_{f_{v,3}}$ and $u_{f_{v,4}}$ going clockwise, with
the edge between faces corresponding to  $u_{f_{v,1}}$ and
$u_{f_{v,2}}$, as well as the edge between faces corresponding to
$u_{f_{v,3}}$ and  $u_{f_{v,4}}$ colored $l_v$, and the other two edges
colored $k_v$.
Then the contribution for the vertex is
\[ \mathrm{Tet}\begin{pmatrix} u_{f_{v,1}} &  u_{f_{v,2}} & k_v \\
 u_{f_{v,3}} &  u_{f_{v,4}} & l_v \end{pmatrix}.\]
This is  well defined because of the symmetry of the
tetrahedral coefficients.

Let $\mathcal{T}(S)$ be the series,

\begin{multline}\label{TS}
\mathcal{T}(S)=\sum_{u_f}(-1)^{\sum_fg_fu_f}
\prod_f \Delta_{u_f}^{x_f}t^{g_fu_f(u_f+2)}\cdot \\
\prod_e\frac{1}{\theta(u_{f_{e,1}},u_{f_{e,2}},k_e)^{v(e)}} 
\prod_v\mathrm{Tet}\begin{pmatrix} u_{f_{v,1}} & u_{f_{v,2}} &
  k_v \\ u_{f_{v,3}} & 
  u_{f_{v,4}} & l_v\end{pmatrix}.
\end{multline}

When $t$ is a $4r$-th root of unity for odd $r$, the sum is finite, and
$\mathcal{T}(S)$ is Turaev's evaluation corresponding to $U_q(sl_2)$.

\begin{theorem} Suppose that $S$ is a shadow in a closed orientable
surface $F$. If the total gleam of $S$ is positive then the series
$\mathcal{T}(S)$  converges when $|t|<1$. If the total gleam
is negative then the series converges when $|t|>1$. The series always
converges when $t=\pm1$ and the genus of $F$ is greater than $1$.
\end{theorem}

\proof 
Let $S$ be a shadow with positive gleam. The case of negative
gleam is similar.
We want to show that for $|t|<1$ the series (\ref{TS})
is absolutely convergent.

The fundamental estimate for tetrahedral coefficients \cite{nefeli} says that
\begin{equation}\label{fundest} 
\left| \text{Tet}\begin{pmatrix} a & b & e \\ c & d & f\end{pmatrix}\right| 
\leq
\sqrt{\frac{\theta(b,c,e)\theta(a,d,e)\theta(a,b,f)\theta(c,d,f)}{(-1)^{e+f}[e+1][f+1]}}.
\end{equation}   
Reorganizing the product according to vertices, the absolute values of
the terms for $\mathcal{T}(S)$ become
\begin{multline}
\prod_f\left| \Delta_{u_f}^{x_f}
t^{g_fu_f(u_f+2)}\right| \cdot \\
\prod_v
\frac{\left|\mathrm{Tet}\begin{pmatrix} u_{f_{v,1}} & u_{f_{v,2}} &
  k_{v} \\ u_{f_{v,3}} & 
  u_{f_{v,4}} & l_{v}\end{pmatrix}\right|}
{\sqrt{\left|\theta(u_{f_{v,1}},u_{f_{v,2}},k_v)
\theta(u_{f_{v,3}},u_{f_{v,4}},k_v)\theta(u_{f_{v,2}},u_{f_{v,3}},l_v) 
\theta(u_{f_{v,4}},u_{f_{v,1}},l_v)\right|}} 
\end{multline}
Applying the estimate (\ref{fundest}), this is less than or equal to
\begin{equation}
\prod_f\left| \Delta_{u_f}^{x_f}t^{g_fu_f(u_f+2)}\right|
\prod_v \frac{1}{\sqrt{[k_v+1][l_v+1]}}.
\end{equation}
Since the product over the vertices does not depend on the coloring of the
faces,
we factor it out to the front as a constant $C$. Now we have
\begin{equation}
C\prod_f\left| \Delta_{u_f}^{x_f}t^{g_fu_f(u_f+2)}\right|.
\end{equation}
Let $B$ (the {\em breadth}) of the shadow be the maximum possible difference of the
colors on any two faces under an admissible coloring. Given a
particular coloring,
let $U$ be the maximal color assigned to a face. Let $p$ be the sum
of the gleams of the faces with positive gleam and $n$ be the sum of the
gleams of the faces with negative gleam, so that $p+n$ is the total gleam of
the shadow. By replacing all $u_f$ on the negative gleam faces by $U$
and all $u_f$ on the faces of positive gleam by $U-B$ we get,
\begin{equation}\label{prod} 
  C\prod_f \left|\Delta_{u_f}^{x_f}t^{g_fu_f(u_f+2)}\right|
\leq 
 C\prod_f \left|\Delta_{u_f}^{x_f}
\right|
t^{nU (U+2)+p(U-B)(U-B+2)}
\end{equation}
We decompose the faces into two subsets, $F^+$ consisting of faces with 
positive Euler characteristic, and $F^-$ consisting of faces with 
negative Euler characteristic. The product (\ref{prod}) is less than
or equal to:
\begin{equation}
  \label{eq:jo}
  C\left|
\frac{
\prod_{f \in F^+}\Delta_{U}^{x_f}}{\prod_{f \in F^-}
\Delta_{U-B}^{-x_f}}
\right||t|^{(n+p)U(U+2)-2pBU+2p(U-B)+pB^2} .
\end{equation}
Rewrite (\ref{eq:jo}) as,
\begin{equation}
  \label{eq:nojo}
  C
\frac{1}{\left|\Delta_{U-B}\right|^{-(2-2g)}}
\prod_{f\in F^+} \left|\frac{\Delta_{U}}{\Delta_{U-B}}\right|^{x_f}
|t|^{(n+p)U(U+2)-2pBU+2p(U-B)+pB^2}
\end{equation}
Since $|t|^{-2n}\leq |\Delta_n|\leq (n+1)|t|^{-2n}$,
\begin{equation}\frac{1}{\left|\Delta_{U-B}^{-(2-2g)}\right|}\leq |t|^{2U(2-2g)},
\end{equation} 
and
\begin{equation}
\prod_{f\in F^+}
\left(\frac{\left|\Delta_{U}\right|}{\left|\Delta_{U-B}\right|}\right)^{x_f}=
\left(\frac{\left|\Delta_{U}\right|}{\left|\Delta_{U-B}\right|}\right)^p\leq 
U^p|t|^{-2Bp}.
\end{equation}
Thus (\ref{eq:nojo}) is less than or equal to
\begin{equation}
CU^p|t|^{2U(2-2g)-2Bp+(n+p)U(U+2)-2pBU+2p(U-B)+pB^2}.
\end{equation}
Since $n+p>0$ the exponent is a quadratic function in $U$ with positive lead
coefficient. Therefore, as $U$ gets large, this is eventually smaller than the
terms of a 
geometric series.

Now reorganize the terms in the sum for $\mathcal{T}(S)$ according to
the maximal color $U$.
There are at most $B^{\#faces}$ admissible colorings
having the same $U$. Thus,
\begin{equation}
  \label{eq:jojo}
  \mathcal{T}(S)\leq \sum_{U}
  B^{\#f}
CU^p|t|^{2U(2-2g)-2Bp+(n+p)U(U+2)-2pBU+2p(U-B)+pB^2 }
\end{equation}
which is clearly convergent.\endproof

\begin{remark}
There exist shadows of total gleam $0$ for which $\mathcal{T}(S)$ is never
convergent when $|t|\neq 1$.
\end{remark}
A simple example is obtained by a separating non null-homotopic closed curve
on a surface of genus $2$ with gleams on the faces equal to $4$ and $-4$.

\section{Limiting behavior at $t=-1$}\label{lim}

In \cite{nefeli} we showed that as $t\rightarrow -1$ from the interior
or the exterior of the unit disk then the Yang-Mills measure converges
to the Yang-Mills measure at $t=-1$. It is the goal of this section to
show that, as $r\rightarrow \infty$, the Yang-Mills measure at the
$4r$-th root of unity, $\ym_r:K_r(F)\rightarrow \mathbb{C}$, converges
to twice the Yang-Mills measure at $t=-1$, 
$\ym_{-1}:K_{-1}(F)\rightarrow \mathbb{C}$  .
Although it can be shown using the Symmetry Principle for the shadow world
formula, we can see this more simply using the formulas from
\cite{nefeli}. Note that in this section, subscripts of $\ym$  are used to
indicate the value of the complex parameter $t$.
 
\begin{remark}
The Yang-Mills measure is
not continuous in the variable $t$ at nontrivial $4r$-th roots of
unity.
\end{remark}
For instance, let $F$ be an annulus and take $s$ to be the skein
determined by $4$ parallel copies of the central core.  Then for
$|t|<1$,  
$\ym_t(s)=2$, while at the $12$-th root of unity $\ym_3(s)=1$.

\begin{theorem}
Let $F$ be a closed surface of genus greater than $1$.
Given a colored framed link $L$ in $F\times I$,  
\begin{equation}\lim_{r\rightarrow\infty}\ym_r(L)=2\cdot \ym_{-1}(L).\end{equation}
\end{theorem}

\proof
The proof requires the use of the following Symmetry Principle for the
Kauffman bracket.

\begin{prop}\label{kauf}
Let $L$ be a framed link, and $K$  a framed knot disjoint from $L$,
with framing $f$. Denote by $K_u$ the knot $K$ colored with $u$, and
by $\langle L\rangle$ the value of the Kauffman bracket of $L$
evaluated at the $4r$-th root of unity.
Then
\begin{equation}\langle L\cup K_u\rangle = i^{f(r+2u-2)}(-1)^{1+\lambda}\langle
L\cup K_{r-u-2}\rangle,\end{equation}
where $\lambda= \sum \mathrm{lk}(L_i,K)$, and the sum is taken over
all components $L_i$ of the link $L$ carrying odd colors.

In particular, if the framing on the component $K$ is $0$ then
\begin{equation}\langle L\cup K_u\rangle =(-1)^{1+\lambda}\langle L\cup K_{r-u-2}\rangle.\end{equation}
\end{prop}

\proof
It follows immediately from the Symmetry Principle proved in
\cite{KM}. Notice that their color $k$ corresponds to our color $k-1$.
\endproof

To compute the Yang-Mills measure of the skein $s$
in $K_r(F)$, first puncture the
surface, to get a surface $F'$ with boundary, and then take the sum over all
the labels $u$ on the 
blackboard framed knot $\partial$ parallel to the boundary \cite{nefeli}, as
follows:
\begin{equation}\label{origin}
\ym_r(s)=\sum_{u=0}^{r-2} \Delta_u\ym_r(\partial_u*s).
\end{equation}
We are using the subscript $r$ in $\ym_r$ to emphasize the fact that
the parameter $t$ is a $4r$-th root of unity.
Here $\ym_r$ on the right hand side is the Yang-Mills measure taken in
the punctured surface $F'$.
Proposition \ref{kauf} implies that (\ref{origin}) is equal to
\begin{equation}\label{half}
2\sum_{u=0}^{\frac{r-3}{2}} \Delta_u\ym_r(\partial_u*s).
\end{equation}
Let $\gamma$ be a trivalent spine of the surface $F'$. Recall that any skein
in $\kr(F')$ can be written as a linear combination of the skeins
corresponding to the $r$-admissible colorings of $\gamma$. Since the
coefficients 
are continuous as $r\rightarrow \infty$, it is enough to show that
\begin{equation}
\lim_{r\rightarrow\infty}\ym_r(s_c)=2\ym_{-1}(s_c),
\end{equation}
where $s_c$ is a skein coming from the coloring $c$ of $\gamma$.

Suppose that $\gamma$ has edges $e$ and vertices $v$. Let $k_e$ denote the
color associated to the edge $e$ in $s_c$, and $k_{v,1}$, $k_{v,2}$,
$k_{v,3}$ be the colors of the three edges ending in a vertex $v$.
Using the formula for $\ym_r(s_c)$ from \cite{nefeli} together with (\ref{half})
we get:
\begin{equation} 
\ym_r(s_c)=2\sum_{u=0}^{\frac{r-3}{2}} \Delta_u \prod_e
\frac{1}{\theta(u,u,k_e)}
\prod_v  \text{Tet}\begin{pmatrix} u & u & u \\ k_{v,1} & k_{v,2} &
k_{v,3}\end{pmatrix}.
\end{equation}
Fix $\epsilon >0$. We will show that there exists $K>0$ so that for all
$r$,
\begin{equation} \label{eq:wojo}
 \left|
   2\sum_{u= K}^{\frac{r-3}{2}} \Delta_u \prod_e
\frac{1}{\theta(u,u,k_e)}
\prod_v  \text{Tet}\begin{pmatrix} u & u & u \\ k_{v,1} & k_{v,2} &
k_{v,3}\end{pmatrix}\right| < \frac{\epsilon}{4},
\end{equation}
and that the analogous inequality, with no upper bound on the summation, holds
for 
$t=-1$. Henceforth we will write the summation without any upper bounds, as the
computation is similar at $t=-1$ and at $t=e^{\frac{\pi i}{2r}}$.
The first step is to rewrite the
summands of (\ref{eq:wojo}) in terms of the vertices, so it becomes,
\begin{equation}
  \label{eq:ohnojo}
 2\sum_{u\geq K} |\Delta_u| \
\prod_v \frac{\left| \text{Tet}\begin{pmatrix} u & u & u \\ k_{v,1} & k_{v,2} &
k_{v,3}\end{pmatrix}\right|}
{\sqrt{\left|\theta(u,u,k_{v,1})\theta(u,u,k_{v,2})
\theta(u,u,k_{v,3})\right|}}.   
\end{equation}
From the fundamental estimate (\ref{fundest}) we see that (\ref{eq:ohnojo}) is
less than or equal to
\begin{equation}\label{almost}
 2\sum_{u\geq K} |\Delta_u| \prod_v \frac{1}{\sqrt{\left|\Delta_u\right|}}
=2\sum_{u\geq K}\frac{1}{[u+1]^{2g-2}},
\end{equation}
since the number of vertices of $\gamma$ is equal to $4g-2$. 
As $[n]\geq \frac{\pi}{2} n$, 
the quantity (\ref{almost}) is less than or
equal to 
\begin{equation}
\left(\frac{2}{\pi}\right)^{2g-2}\sum_{u\geq K}\frac{1}{(n+1)^{2g-2}},
\end{equation}
which is smaller than $\frac{\epsilon}{4}$ for big enough $K$.

Given $K$,
\begin{equation}
\left|2\sum_{u=0}^K \Delta_u\ym_r(\partial_u*s) -
2\sum_{u=0}^K \Delta_u\ym_{-1}(\partial_u*s)
\right|
\leq\frac{\epsilon}{2}
\end{equation}
if $r$ is big enough.
Thus
\begin{multline}
\left|\ym_r(s_c)-2\ym_{-1}(s_c)\right|\\
\leq \left|2\sum_{u=0}^K \Delta_u\ym_r(\partial_u*s) -
2\sum_{u=0}^K \Delta_u\ym_{-1}(\partial_u*s)
\right| \\
+ \left| 2\sum_{u= K}^{\frac{r-3}{2}} \Delta_u \prod_e
\frac{1}{\theta(u,u,k_e)}
\prod_v  \text{Tet}\begin{pmatrix} u & u & u \\ k_{v,1} & k_{v,2} &
k_{v,3}\end{pmatrix} \right| \\
+ \left|2\sum_{u=K}^{\infty} \Delta_u\ym_{-1}(\partial_u*s)   \right|\leq
\epsilon. 
\end{multline}
\endproof

\section{Integrating on the character variety}

We will use the formula from Theorem \ref{theformula} to compute some
examples. The irreducible $SU(2)$-characters of
$\pi_1(F)$ will be denoted $\mathrm{char}^i(\pi_1(F))$.
Recall that the algebra $K_{-1}(F)$ is
isomorphic to 
the ring of polynomial functions on  $\mathrm{char}^i(\pi_1(F))$,
\cite{bu, PS}. The Yang-Mills measure is the linear
functional induced by integrating against the symplectic measure \cite{nefeli}.
Given a diagram on $F$ it corresponds to a function
on $\mathrm{char}^i(\pi_1(F))$.
Under this correspondence a knot on a surface is taken to the function
that assigns to a representation $\rho$ the negative of the trace of a loop
determined by the knot in $\pi_1(F)$.
In addition, 
any colored diagram represents an element of $H_1(F;{\mathbb Z}_2)$ by
replacing each edge colored with $n$ by $n$ parallel strands, hooked
together at trivalent vertices as if they were triads.
\begin{prop}
If the diagram $\alpha$ does not represent $0$ in  $H_1(F;{\mathbb
  Z}_2)$ then integrating  $\alpha$ over $\mathrm{char}^i(\pi_1(F))$
yields zero. In particular, if $\alpha$ is a non-separating curve,
then the integral of $\alpha$ is $0$.
\end{prop}
\proof
The diagram of $\alpha$ has no admissible colorings. 
If $\alpha$ represents a non-zero element of $H_1(F;{\mathbb Z}_2)$
then there is a path $\gamma$ in $F$ which has non-zero intersection
number with the cycle represented by $\alpha$. In an admissible
coloring, the difference of the colors on two adjacent faces is congruent
modulo $2$ to the color of the edge between them. Hence the
intersection number of $\alpha$ and $\gamma$ is the sum of the
differences of the colors of the faces that $\gamma$ passes through. This
contradicts the fact that $\gamma$ has a nonzero 
${\mathbb Z}_2$-intersection number with $\alpha$. \endproof

\begin{prop}
Let $\beta$ be a simple closed curve which separates the closed surface
$F$ into surfaces of genus $G_1>0$ and $G_2>0$ respectively. Integrating
 $\beta$ yields
\begin{equation}\ym(\beta)
=-\sum_{u=1}^{\infty}\left(\frac{1}{u^{2G_1-1}(u+1)^{2G_2-1}} +
\frac{1}{(u+1)^{2G_1-1}u^{2G_2-1}}\right).
\end{equation}
\end{prop}
\proof
The diagram of $\beta$ has two faces, $f_1$ and
$f_2$, meeting along the only edge $e$  colored with $1$. There are no
vertices, so $v(e)=0$.
Let $u_i$ be the color assigned to $f_i$, for $i=1,2$. In order for
the triple $(u_1, u_2,1)$ to be admissible, we must have  $u_1=u_2-1$
or $u_1=u_2+1$. Thus we  have a sum over all nonnegative integers
$u$ with two 
terms: one where $u_1=u$ and $u_2=u+1$, and one where $u_1=u+1$ and
$u_2=u$. The Euler characteristics of the faces are
$x_1=1-2G_1$ and $x_2=1-2G_2$. The gleams  $g_1$ and $g_2$ of both
faces are  zero.

When $t=-1$, all the quantized integers in our formulas are replaced
by ordinary integers.  We have:

\begin{multline}\ym(\beta)= \\\sum_u \left(
\Delta_u^{1-2G_1}
\Delta_{u+1}^{1-2G_2}
\frac{1}{(\theta(u,u+1,1))^0}
+\Delta_{u+1}^{1-2G_1}
\Delta_{u}^{1-2G_2}
\frac{1}{(\theta(u+1,u,1))^0}\right)
\end{multline} 
Note that $\Delta_u=(-1)^u(u+1)$. 

Substituting we get
\begin{multline}\ym(\beta)=\\ \sum_u
\left( (-1)^{p_1}\frac{1}{(u+1)^{2G_1-1}(u+2)^{2G_2-1}}
+(-1)^{p_2}\frac{1}{(u+2)^{2G_1-1}(u+1)^{2G_2-1}}\right) \\
=-\sum_{u=1}^{\infty}\left(\frac{1}{u^{2G_1-1}(u+1)^{2G_2-1}} +
\frac{1}{(u+1)^{2G_1-1}u^{2G_2-1}}\right),\end{multline}
since
$p_1=u(1-2G_1)+(u+1)(1-2G_2)$
and $p_2=(u+1)(1-2G_1)+u(1-2G_2)$ are both odd.\endproof

We can express $\ym(\beta)$ in terms of Bernoulli numbers.

\begin{cor}
Suppose that $G_1=G_2=G$. Then
\begin{equation}
\ym(\beta)=
-2\sum_{j=1}^{2G-1}\binom{4G-j-3}{2G-j-1}
 +4\sum_{j=1}^{G-1}
\binom{4G-2j-3}{2G-2j-1} 
 \frac{1}{2}\left|B_{2j}\right|\frac{(2\pi)^{2j}}{(2j)!}
\end{equation}
Otherwise let $G=\text{min}(G_1,G_2)$, and $H=\text{max}(G_1,G_2)$.
Then  
\begin{equation}
\ym(\beta)=
-\sum_{j=1}^{2G-1}\binom{2G+2H-j-3}{2G-j-1}
-\sum_{j=1}^{2H-1}\binom{2G+2H-j-3}{2H-j-1}+
\end{equation}
\[2\sum_{i=1}^{G-1}
\left(\binom{2G+2H-2i-3}{2G-2i-1}+\binom{2G+2H-2i-3}{2H-2i-1}\right)
\frac{1}{2}\left|B_{2i}\right|\frac{(2\pi)^{2i}}{(2i)!}\]
\[+2\sum_{i=G}^{H-1}\binom{2G+2H-2i-3}{2H-2i-1}     
\frac{1}{2}\left|B_{2i}\right|\frac{(2\pi)^{2i}}{(2i)!}.
\]
\end{cor}
\proof
Collect the terms and note that if $r$ is even, then
$\sum\frac{1}{n^r}=\frac{1}{2}\left|B_r\right|\frac{(2\pi)^r}{r!}$,
where $B_r$ is a Bernoulli number.
\endproof

We finish by giving numerical values in Figure \ref{integrals} for the
Yang-Mills measure of 
separating curves on surfaces with genus up to $10$.
Here $G=G_1+G_2$ is the genus of the surface, the curve splits it into
surfaces of genus $G_1$ and $G_2$ respectively.
\begin{figure}
\newcommand{\fu}{{\vrule height17pt width0pt depth11pt}}
\begin{center}
\begin{tabular}{|l|l|l|l|}
\hline
G & $G_1$ & $G_2$ & Value of the integral \\ \hline
$2$\fu & 1 & 1 & $-2$ \\ \hline
$3$\fu & 1 & 2 & $ -4+{\displaystyl\frac{1}{3}}\pi^{2}$ \\ \hline
$4$\fu & 1 & 3 &
$-6+{\displaystyl\frac{1}{45}}\pi^{4}+
{\displaystyl\frac{1}{3}}\pi^{2}$ \\ \hline
$5$\fu & 1 & 4 & $-8+{\displaystyl\frac{2}{945}}\pi^{6} + 
{\displaystyl\frac{1}{45}}\pi^{4} +
{\displaystyl\frac{1}{3}}\pi^{2}$ \\ \hline
$6$\fu & 1 & 5 & $ -10 + {\displaystyl\frac{1}{4725}}\pi^{8}
+ {\displaystyl\frac{2}{945}}\pi^{6} +
{\displaystyl\frac{1}{45}}\pi^{4} +
{\displaystyl\frac{1}{3}}\pi^{2}$ \\ \hline
$4$\fu & 2 & 2 & $ -20+2\pi^{2}$ \\ \hline
$5$\fu & 2 & 3 &
$-56+{\displaystyl\frac{1}{15}}\pi^{4}+5\pi^{2}$ \\ \hline
$6$\fu & 2 & 4 & $ -120+{\displaystyl\frac{2}{315}}\pi^{6} + 
{\displaystyl\frac{2}{9}}\pi^{4} +
{\displaystyl\frac{28}{3}}\pi^{2}$ \\ \hline
$7$\fu & $2$ & $5$ & $ -220
+{\displaystyl\frac{1}{1575}}\pi^{8} +   
{\displaystyl\frac{4}{189}}\pi ^{6} +
{\displaystyl\frac{7}{15}}\pi^{4} + 15\pi^{2}$ \\ \hline
$6$\fu & 3 & 3 & $- 252+{\displaystyl\frac{2}{9}}\pi^{4} + 
{\displaystyl\frac{70}{3}}\pi^{2}$ \\ \hline
$7$\fu & $3$ & $4$ &  $-792+{\displaystyl\frac{2}{189}}\pi^{6} + 
{\displaystyl\frac{14}{15}}\pi^{4} + 70\pi^{2}$\\ \hline
$8$\fu & $3$ & $5$ & $ - 2002+{\displaystyl\frac{1}{945}}\pi^{8}
+  {\displaystyl\frac{2}{27}}\pi^{6} + 3\pi^{4}+165\pi^{2}$ \\ \hline
$8$\fu & $4$ & $4$ & $- 3432
+{\displaystyl\frac{4}{135}}\pi^{6} +  
{\displaystyl \frac {56}{15}}\pi^{4} + 308\pi^{2}$ \\ \hline
$9$\fu & $4$ & $5$ & $ -
11440+{\displaystyl\frac{1}{675}}\pi^{8} + 
{\displaystyl\frac{62}{315}}\pi^{6} +
{\displaystyl\frac{209}{15}}\pi^{4} + 1001\pi^{2}$ \\ \hline
$10$\fu & $5$ & $5$ & $ - 48620
+{\displaystyl\frac{2}{525}}\pi^{8}  + 
{\displaystyl\frac{44}{63}}\pi^{6} +
{\displaystyl\frac{286}{5}}\pi^{4} + 4290\pi^{2}$ \\ \hline
\end{tabular} 
\end{center}
\caption{Integrals of separating curves on a closed surface of
  genus $G$ }\label{integrals}
\end{figure}

\Addresses

\end{document}